\documentclass[letterpaper, 11pt, conference]{ieeeconf}  % Comment this line out
\IEEEoverridecommandlockouts                              % This command is only
\usepackage{amssymb, graphicx, url}
\def\R{\mathbb{R}}
\newtheorem{theorem}{Theorem}

\title{\LARGE \bf
Fast Algorithms for Distributed Optimization and Hypothesis Testing: A Tutorial
}

\author{Alex Olshevsky% <-this % stops a space
\thanks{ This research was
supported by NSF under award CMMI-1463262 and AFOSR under award FA-95501510394}% <-this % stops a space
\thanks{A. Olshevsky is with the department of Electrical and Computer Engineering at Boston University.  {\tt\small alexols@bu.edu}}%
\thanks{This is a corrected version of an invited paper that appeared in the Proceedings of CDC 2016.}
}

\begin{document}

\maketitle
\thispagestyle{empty}
\pagestyle{empty}

%%%%%%%%%%%%%%%%%%%%%%%%%%%%%%%%%%%%%%%%%%%%%%%%%%%%%%%%%%%%%%%%%%%%%%%%%%%%%%%%
\begin{abstract} We consider several problems in the field of distributed optimization and  hypothesis testing. We show how to obtain  convergence times for these problems that scale linearly with the total number of nodes in the network by using a recent linear-time algorithm for the  average consensus problem. 
\end{abstract}

%%%%%%%%%%%%%%%%%%%%%%%%%%%%%%%%%%%%%%%%%%%%%%%%%%%%%%%%%%%%%%%%%%%%%%%%%%%%%%%%
\section{Introducion}

 Throughout we will assume we are given a network of nodes, modeled by a connected undirected graph $G= (\{1, \ldots, n\}, \mathcal{E})$, and that nodes can send messages to their neighbors in this graph.  The goal of the network is to collectively achieve some global coordination task based only on local interactions. We will focus on three simple such tasks.

\bigskip

\noindent {\bf Fusion of Gaussian Measurements:} Node $i$ has access to the noisy measurement \[ y_i  = \theta + w_i \] of some parameter of interest $\theta \in \R$. Here each $w_i$ is a $N(0, \sigma_i^2)$ Gaussian noise. The nodes desire to compute a maximum likelihood estimate of the unknown parameter $\theta$.

\bigskip

\noindent {\bf Distributed Separable Optimization:} Each node $i$ knows a convex function $f_i(\cdot): \R \rightarrow \R$ and the nodes desire to collectively agree on a minimizer of 
\[ F(z) = \sum_{i=1}^n f_i(z). \] The problem was first proposed in \cite{angas}.

\bigskip

\noindent {\bf Distributed Learning:} There is a finite set $\Theta = \{\theta_1, \ldots, \theta_k\}$ which describes $k$ possible ``states of the world.'' Each node receives an infinite stream of i.i.d. measurements, and these measurements would be differently distributed in different states of the world.

Specifically, node $i$ receives the infinite sequnce  of random variables \[ S_1^i, S_2^i, \ldots \] If $\theta_p \in \Theta$ was the true state of the world, these random variables would be i.i.d. with density which we denote by $\ell^i(\cdot ~|~ \theta_p)$.

The goal of every node is to figure out which state of the world best represents {\em all} the samples in the network. There may be a single true state of the world, in which case the nodes would like to figure out what it is; if not, the nodes would like to figure out which state(s) of the world best model all the data the network keeps collecting.

More specifically, each node would like to maintain a ``belief vector'' which assigns a number to every hypothesis in $\Theta$; over time, as the system collects more and more measurements, this belief vector should concentrate on the best hypotheses.  This problem was first proposed in \cite{firstlearn}.

\bigskip

All of the above problems are closely related to the so-called average consensus problem. In the rest of this tutorial, we will introduce the average consensus problem, discuss the best available algorithms for it, and discuss the implications of these algorithms for the above three problems.

\section{The average consensus problem and the associated convergence times}

Given that node $i$ in an undirected graph begins with a real number $x_i(0) \in \R$, the average consensus protocol asks for a distributed algorithm by means of which the nodes can all agree on the average $(1/n) \sum_{i=1}^n x_i(0)$. It is assumed that time is broken up into a discrete sequence of rounds $t=1, 2, 3, \ldots$ and that, during each round, nodes can send each other real-valued messages as well as perform simple computations.

Of course, in reality nodes will not exchange real numbers but rather quantized versions of these numbers. Because of this it is important that protocols for average consensus perform only computations that behave well under perturbation and truncation. Indeed, all the protocols we discuss in this tutorial are based on simple linear updates.

The simplest idea is for each node to iterate as 
\[ x_i(t+1) = \sum_{j \in N(i)} a_{ij} x_j(t), \] where $N(i)$ denotes the set of neighbors of agent $i$ in the graph. We will adopt the convention throughout this paper that $i \in N(i)$, i.e., every node is its own neighbor.

Under appropriate conditions, it can be shown that this works. Indeed, for convenience, we adopt the notation 
\[ \overline{x} =  \frac{\sum_{j=1}^n x_j(0)}{n}. \] We then have the following theorem.

\bigskip

\begin{theorem} \label{conv} Suppose the matrix $A = [a_{ij}]$ satisfies the following properties:
\begin{itemize} \item $A$ is doubly stochastic (meaning that every entry of $A$ is nonnegative and every row and every column adds up to $1$).
\item If $(i,j) \in \mathcal{E}$, then $a_{ij} > 0$ and $a_{ji} > 0$. 
\end{itemize} Then for every $i = 1, \ldots, n$, we have 
\[ \lim_{t \rightarrow \infty} x_i(t) = \overline{x}. \]
\end{theorem}

\bigskip

\bigskip

Note that the second bullet implies that the diagonal of $A$ is strictly positive, as by our convention $(i,i) \in \mathcal{E}$ for all $i$.

The theorem in this form is taken from \cite{othesis} but is not origial to that thesis; rather, it is a minor variation on the earlier results of \cite{degroot, tsi, JLM} .

There are a number of ways to choose the matrix $A$ satisfying the assumptions of the theorem. The easiest is for each node $i$ to set $a_{ij} = \epsilon$ for all $j \in N(i)$, and then set $a_{ii} = 1 - d(i) \epsilon$, where $d(i)$ is the degree of node $i$. This works when $\epsilon$ is strictly smaller than the largest degree in the network.

A corresponding convergence time bound has been obtained in \cite{four}. Let us stack up the numbers $x_1(t), \ldots, x_n(t)$ into the vector ${\bf x}(t)$, and we will use ${\bf 1}$ to denote the all-ones vector. Finally, let us adopt the notation 
\[ E(t) = \left| \left| {\bf x}(t) - \overline{x} {\bf 1} \right| \right|_2^2. \] The quantity $E(t)$ measures how far we are from average consensus; we will refer to it as the {\em squared error} at time $t$. We then have the following theorem from \cite{four}.

\bigskip

\bigskip

\begin{theorem} Under the same assumptions as Theorem \ref{conv}, and additionally adopting the notation 
\[ \eta = \min_{a_{ij} > 0} a_{ij} \] we have that 
\[ E(t) \leq \left( 1 - \frac{\eta}{2n^2} \right)^t E(0).  \]
\end{theorem}

\bigskip

\bigskip

In other words, the squared error contracts by a factor of $1-\eta/(2n^2)$ at each step. This implies that every
$2n^2/\eta$ steps, $E(t)$ shrinks by a factor of $e^{-1}$. In other words, it takes $(2n^2/\eta) \log (1/\epsilon)$ steps for $E(t)$ to shrink to $\epsilon$ times its initial value.

It may be remarked that this bound depends on $\eta$, and in a sense something like this is unavoidable; if all the off-diagonal entries of $A$ are close to zero, convergence should be slow.

Referring back to the construction where every node places a value of $\epsilon$ on its neighbors, we have that if $\epsilon = 1/2 d_{\rm max}$ (where $d_{\rm max}$ is the largest degree, the last theorem implies a time of 
\[ O \left( n^2 d_{\rm max} \log \frac{1}{\epsilon} \right) \] until the squared error $E(t)$ shrinks by a factor of $\epsilon$.

\bigskip

One way to obtain a nicer bound is to rely on the so-called lazy Metropolis weights. Each node $i$ sets 
\[ a_{ij} = \frac{1}{2 \max(d(i), d(j))},  \mbox{ for all } j \in N(i), ~~ j \neq i \] and 
\[ a_{ii} = 1 - \sum_{j \in N(i), ~~ j \neq i} a_{ij}. \] We then have the folowing theorem, which follows by some elementary manipulations from \cite{linear}.

\bigskip

\bigskip

\begin{theorem} If $a_{ij}$ are set to be the lazy Metropolis weights then 
\[ \left| \left| {\bf x}(t) - \overline{x} {\bf 1} \right| \right|_2^2 \leq \left( 1 - \frac{1}{37 n^2} \right)^t \left| \left| {\bf x}(0) - \overline{x} {\bf 1} \right| \right|_2^2  \]
\end{theorem}

\bigskip

\bigskip

Using this theorem, the time until the squared error shrinks by a factor of $\epsilon$ is 
\[ O \left( n^2 \log \frac{1}{\epsilon} \right) \]
with the lazy Metropolis weights.

One can improve this further to a {\em linear} scaling in the number of nodes, provided that a reasonably accurate upper bound on $n$ (the total number of nodes in the network) is known to all the nodes.

{\em Indeed, suppose every node knows a number $U$ such that $U \geq n$.} Then nodes can implement the following update.

\begin{eqnarray} y_i(t+1) & = & x_i(t) + \sum_{j \in N(i)} \frac{x_j(t) - x_i(t)}{2 \max(d(i), d(j))} \nonumber \\ 
x_i(t+1) & = & y_i(t+1) \label{am} \\ && + \left( 1 - \frac{2}{9U+1} \right) \left( y_i(t+1) - y_i(t) \right) \nonumber
\end{eqnarray}

We then have the following theorem, which again follows by some elementary manipulations from \cite{linear}.

\bigskip

\begin{theorem} Suppose the nodes implement the iteration of Eq. (\ref{am}). Then 
\[ E(t) \leq 18 \left( 1 - \frac{1}{9U} \right)^t E(0). \]
\end{theorem}

\bigskip

\bigskip

Under this theorem, the time until $E(t)$ shrinks by a factor of $\epsilon$ is 
\[ O \left( U \log \frac{1}{\epsilon} \right) \] under the protocol of Eq. (\ref{am}). If $U=n$ or $U$ is within a constant factor of $n$, this convergence time is linear in the number of nodes.

\section{Fast Algorithms for Distributed Optimization and Hypothesis Testing}

We now revisit each of the problems we started with and show how they may be solved using algorithms for average consensus.

\bigskip

\bigskip

\noindent {\bf Fusion of Gaussian Measurements:} In this case, the reduction is immediate; one can check that the maximum likelihood estimate $\widehat{\theta}$ is 
\[ \widehat{\theta} = \frac{\sum_{i=1}^n y_i/\sigma_i^2}{\sum_{i=1}^n 1/\sigma_i^2} = \frac{(1/n) \sum_{i=1}^n y_i/\sigma_i^2}{(1/n) \sum_{i=1}^n 1/\sigma_i^2}. \] Using the linear-time consensus scheme described in the last section, we can compute arbitrarily accurate approximations to the numerator and denominator in linear time.

\bigskip

\bigskip

\noindent {\bf Distributed Optimization:} Our starting point is the algorithm of \cite{angas}, which had each node update as 
\begin{equation} \label{AngAsu} x_i(t+1) = \sum_{j \in N(i)} a_{ij} x_j(t) - \alpha s_i(t), \end{equation} where, as before, $A = [a_{ij}]$ has to be a stochastic matrix and $s_i(t)$ is the subgradient of $f_i(\cdot)$ at $x_i(t)$.

Approximation guarantees which go to zero as $\alpha \rightarrow 0$ (under appropriate technical assumptions such as boundedness of subgradients) were shown in \cite{angas}.

Intuitively speaking, what the nodes really want to do is to do gradient descent on the average function $(1/n) \sum_{i=1}^n f_i(\cdot)$. Unfortunately, every node only knows its own function, making this impossible. The solution is to ``run through'' all the locally-obtained gradients through an average consensus scheme, and this is essentially what Eq. (\ref{AngAsu}) does.

This intuition suggests we may obtain improved performance by using the linear time consensus scheme, and indeed this turns out to be the case as we now discuss. Based on this idea, \cite{linear} proposed the following scheme 
\begin{scriptsize}
\begin{eqnarray} y_i(t+1) & = & x_i(t) +  \frac{1}{2} \sum_{j \in N(i)} \frac{x_j(t) - x_i(t)}{\max(d(i), d(j))} - \beta g_i(t)  \nonumber \\ 
z_i(t+1) & = & y_i(t) - \beta g_i(t) \label{optaccel} \\ 
x_i(t+1) & = & y_i(t+1) + \sigma \left( y_i(t+1) - z_i(t+1) \right) \nonumber
\end{eqnarray} \end{scriptsize} where $\sigma = (  1 - {2}/({9U+1}) )$, $g_i(t)$ is the subgradient of $f_i(\theta)$ at $\theta = y_i(t)$ and $\beta$ is a step-size to be chosen later.

Intuitively speaking, what this scheme does is replace the average consensus scheme coming from a doubly stochastic matrix with the improved, linear time consensus scheme.

We now describe the available results from \cite{linear} on this scheme; the following theorem and the definitions preceeding it are reproduced from \cite{linear}. We will use ${\cal W^*}$ to denote the set of global minima of $f(\theta)$. We will use the standard hat-notation for a running average, e.g., 
 $\widehat{y}_i(t) = (1/t) \sum_{k=1}^t y_i(k)$. To measure the convergence speed of our protocol, we introduce two measures of performance. One is the dispersion of a set of points, intuitively measuring how far from each other the points are, \[ {\rm Disp}(\theta_1, \ldots, \theta_n) = \frac{1}{n} \sum_{i=1}^n \left| \theta _i - {\rm median}(\theta_1, \ldots, \theta_n) \right| \] The other is the ``error'' corresponding to the function $f(\theta) = (1/n) \sum_{i=1}^n f_i(\theta)$, 
\[ {\rm Err}(\theta_1, \ldots, \theta_n) = \left( \frac{1}{n} \sum_{i=1}^n f_i (\theta_i) \right) - f(w^*) \] where $w^*$ is any point in the optimal set ${\cal W}^*$. Note that to define ${\rm Err}(\theta_1, \ldots, \theta_n)$ we need to assume that the set of global minima ${\cal W}^*$ is nonempty. We then have the following theorem from \cite{linear}.

\bigskip

\begin{theorem} Suppose that $U \geq n$ and $U=\Theta(n)$, that $w^* \in {\cal W}^*$, and that the absolute value of all the subgradients of all $f_i(\theta)$ is bounded by some constant $L$. If every node implements the update of Eq. (\ref{optaccel}) with $\beta = \frac{1}{L \sqrt{UT}}$ we then have
\begin{footnotesize} \begin{eqnarray*}
{\rm Disp}(\widehat{y}_1(T), \ldots, \widehat{y}_n(T)) & = & O \left(  \sqrt{\frac{n}{T}} \left(  \frac{||y(0) - \overline{x} {\bf 1}||_2}{\sqrt{T}} + 1 \right) \right)   \\
 {\rm Err}(\widehat{y}_1(T), \ldots, \widehat{y}_n(T)) & = & O \left( L \sqrt{\frac{n}{T}}  \left(  \frac{||y(0) - \overline{x} {\bf 1}||_2}{\sqrt{T}} + 1 \right. \right. \\ && \left. \left.  + (\overline{x} - w^*)^2 \right) \right)  
\end{eqnarray*} \label{opthm}
\end{footnotesize}

\end{theorem}

Focusing on scaling with $n$ and $\epsilon$, we remark that by taking $T=O(n/\epsilon^2)$ we can make the upper bounds on the right-hand side be below $\epsilon$. That is to say, convergence time for this scheme is linear in the number of nodes.

\bigskip

\bigskip

\noindent {\bf Decentralized Hypothesis Testing.} We begin with a discussion of a nearly identical class of update rules for this problem that were proposed in \cite{ucsd, rakh, us}. Each node 
$i$ maintains a stochastic vector $\mu_i^t$ whose cardinality is $|\Theta|$. We will refer to $\mu_i^t$ as the belief of node $i$ and time $t$.  We will use $\mu_i^t (\theta_k)$ to refer to the $k$'th element of $\mu_i^t$. Informally speaking, we mean $\mu_i^t(\theta_k)$ to capture the probability that agent $i$ assigns to $\theta_k$ being the correct hypothesis. We will assume that initially, all agents believe all hypotheses are equally likely, i.e., $\mu_i^0(\theta_k) = 1/|\Theta|$ for all $\theta_k \in \Theta$.

We do not assume that any of the hypotheses in $\Theta$ represent the true state of the world. Our goal will be for the agents to have beliefs which concentrate on the best hypotheses. Following the convention in the hypothesis testing literature, ``best'' is defined in terms of KL divergence. Specifically, let $f^i(\cdot)$ be the true distribution of the measurements of agent $i$ (i.e., $S_i^1, S_i^2, \ldots$ are all distributed acoording to $f^i$). For each $\theta \in \Theta$, we define 
\[ C(\theta) = \sum_{i=1}^n D_{\rm KL} (f^i ~||~ \ell^i(\cdot ~|~ \theta) ). \] Recall that $\ell^i(\cdot ~|~ \theta)$ would be the distribution of measurements received by agent $i$ if $\theta$ were the correct state of the world. Thus our goal is for all agents agents to concentrate their beliefs on those $\theta \in \Theta$ that minimize $C(\theta)$.

We will adopt the convention of writing

\[ \mu_i^{t+1}(\theta) \sim q( \mu_i^t(\theta)), ~~~~ \mbox{ for all } \theta \in \Theta,  \] where $q$ is some function that maps nonnegative numbers to nonnegative numbers to mean that $\mu_i^{t+1}$ is obtained by normalizing $q(\mu_i^t)$ so that its entries add up to one. Adding normalizations will make our formulas rather messy, which is why it is neater to sweep them under the rug with the $\sim$ symbol.

The papers \cite{rakh, ucsd, us} studied variations on the dynamics

\begin{figure*}[ht] 
\begin{tabular}{cc} 
\includegraphics[scale=0.45]{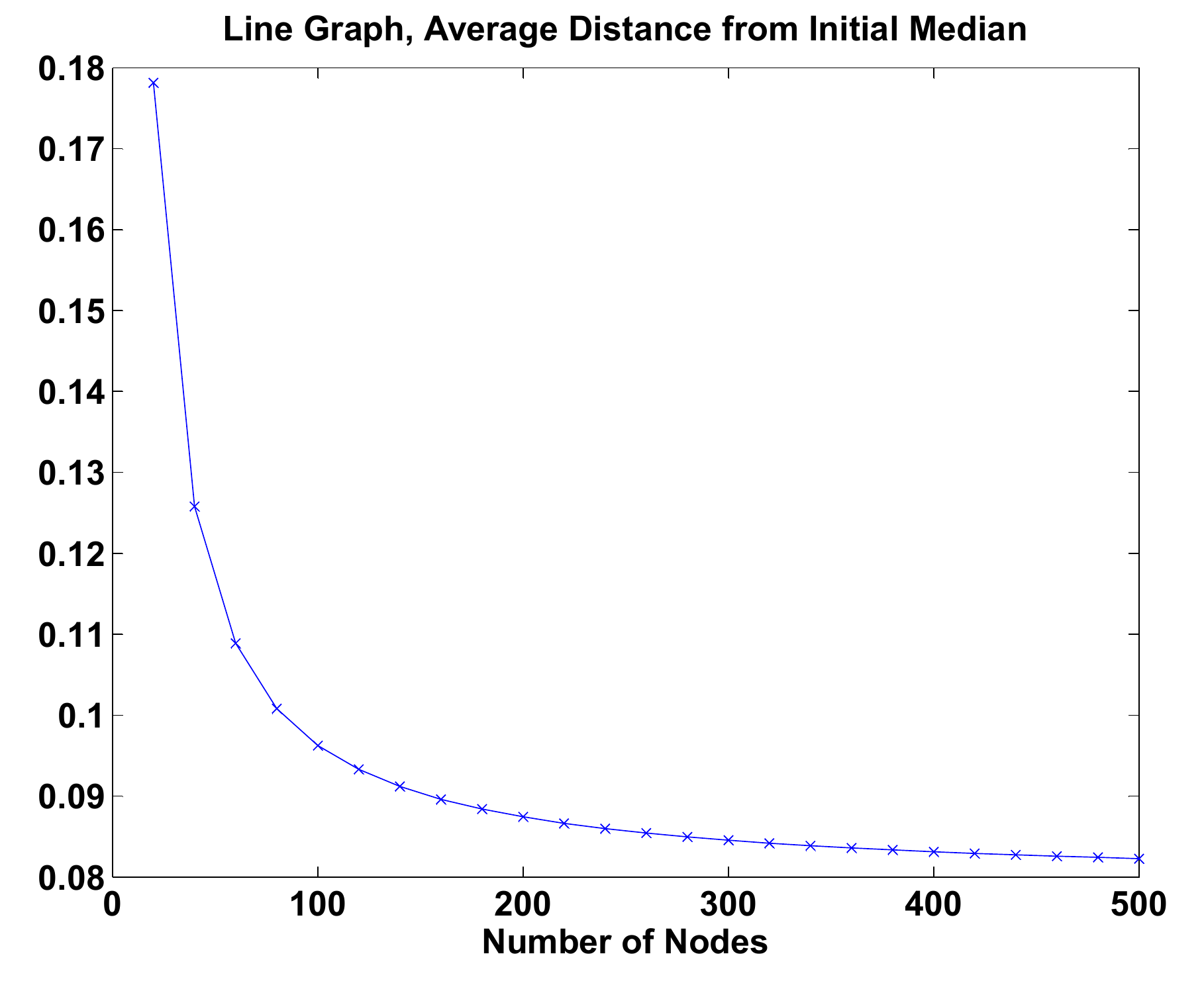} 
&  \includegraphics[scale=0.45]{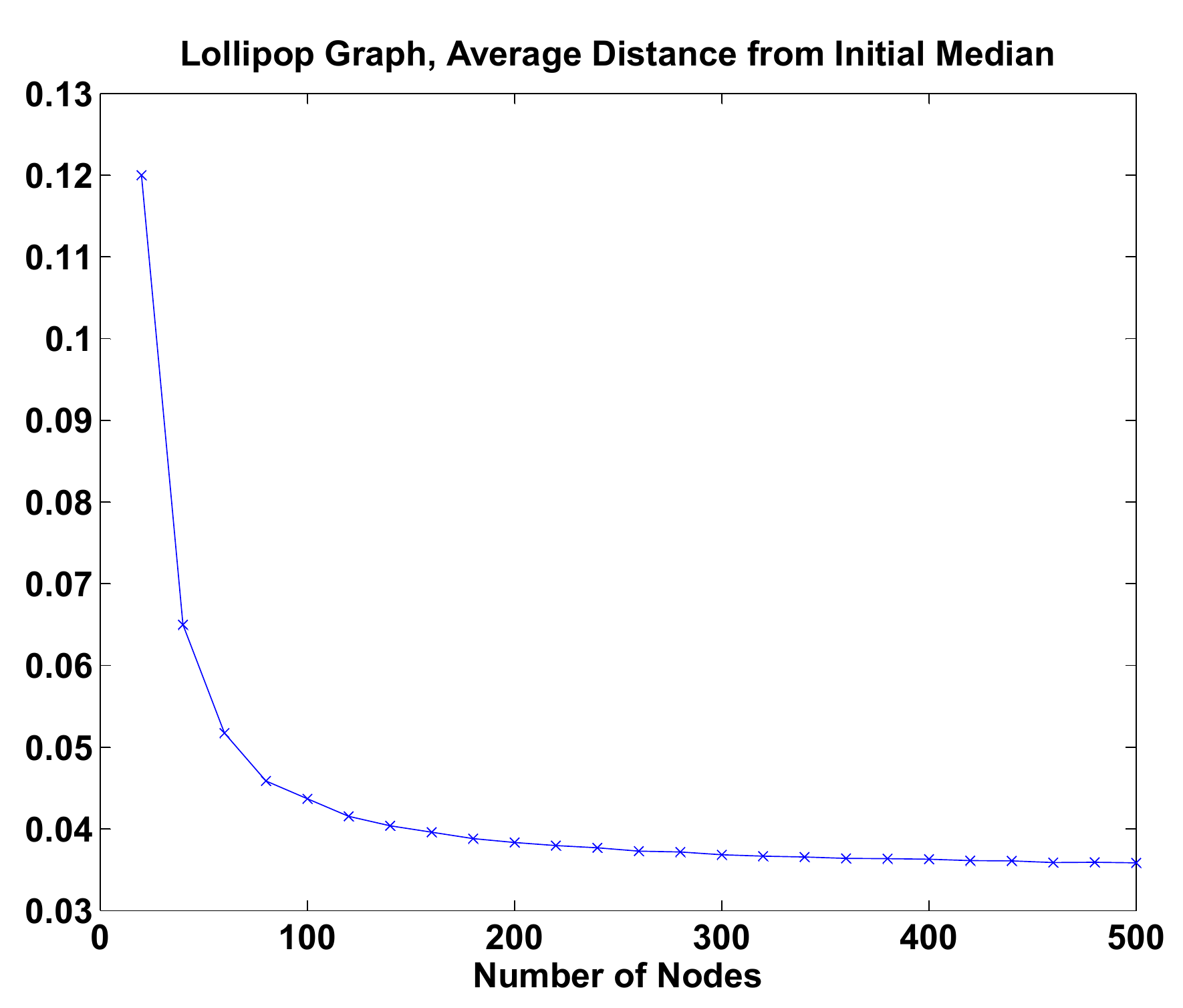} 
\end{tabular} \caption{Average deviation from the median as a function of the number of nodes using the decentralized optimization protocol of Eq. (\ref{optaccel}) after $T=4n$ iterations. The line graph is shown on the left and and convergence time the lollipop graph is shown on the right. The figures are reproduced from \cite{linear}.} \label{opt-simul}
\end{figure*}

\begin{equation} \label{learndyn} \mu_i^{t+1}(\theta) \sim   \ell^i(S_i^{t+1} ~|~ \theta) \prod_{j \in N(i)} \mu_j^t(\theta)^{a_{ij}}, \end{equation} where, as before, $A = [a_{ij}]$ is a doubly stochastic matrix. It was shown in those papers that, under appropriate technical assumption, this scheme works; namely that  the beliefs of each agent concentrate on the correct hypothesis if there is one, or on the set of best hypotheses otherwise.

These dynamics are  inspired by Bayes rule. Essentially, every node performs a Bayes-like update, treating the beliefs of its neighbors as observations (after raising them to the power $a_{ij}$), and making no effort to account for the dependency among the messages it receives.

It is possible to connect these dynamics to the consensus iteration. Indeed, taking logs of both sides of Eq. (\ref{learndyn}), we observe that the quantities $\log \mu_i^t$ are updated via consensus iteration which is ``driven'' by the ``exogenous inputs'' $\log \ell^i(S_i^{t+1} ~|~ \theta)$.

Keeping to this intuition, we may obtain an improvement using our accelerated consensus scheme. Thus \cite{us} proposed the update
\begin{equation} \label{betu} \mu_i^{t+1}(\theta) \sim \frac{\ell^i (S_i^{t+1} ~|~ \theta) \prod_{j \in N(i)} \mu_j^t(\theta)^{(1+\sigma) b_{ij} }}{\prod_{j \in N(i)} \left( \mu_j^{t-1}(\theta) \ell_j^{t} (S_j^t ~|~ \theta) \right)^{\sigma b_{ij}}} \end{equation} where $B = [b_{ij}]$ is the matrix of Metropolis weights and $\sigma = 1 - 2/(9U+1)$. We also initialize $\mu_i^{-1} = \mu_i^0$. Although this scheme appears quite involved, after taking logs of both sides and performing some elementary manipulations, one may verify that what it does is simply replace the ``typical'' doubly-stochastic consensus scheme with the more involved linear-time scheme from \cite{linear}.

We next describe a convergence time guarantee for this scheme. Let us introduce the notation $g$ to denote the gap between the smallest $C(\theta)$ and the second smallest $C(\theta)$ as $\theta$ runs over the finite set $\Theta$.  Moreover, let us make the technical assumption that if $f^i(s) > 0$ then $\ell^i(s ~|~ \theta) > \alpha$ for all $\theta$. We then have the following theorem from \cite{us}.

\bigskip

\bigskip

\begin{theorem} Suppose $U \geq n$ and that $\mu_i^t$ are updated according to Eq. (\ref{betu}). Given a number $\rho$, let us define 
\[ N(\rho) = \lceil \frac{48 (\ln \alpha)^2 \ln (1/\rho)}{(g/n)^2} \rceil. \] Then, for any $\rho \in (0,1)$ the following statement is true: with probability $1-\rho$ we have that for all $t \geq N(\rho)$ and all $\theta_v \not \in \arg \min_{\theta} C(\theta)$, 
\[ \mu_i^t(\theta_v) \leq e^{-(g/n) (t/2) + Z} \] 
where 
\[ Z = O \left( U (\log n) \left( \log \frac{1}{\alpha} \right) \right) \]
\end{theorem}

\bigskip

\bigskip

We may think of $g/n$ as the average ``learning rate'' in the network. Indeed, recall that $C(\theta)$ is a sum of $n$ terms, and $g$ is the gap between the smallest and second-smallest $C(\theta)$; thus $g/n$ may be intuitively thought as the average rate at which the network closes the gap between the first-best and second-best hypotheses. The above theorem tells us that $(g/n)t$ has to be larger than a transient of whose size is $O(U \log n)$ (taking $\alpha$ to be a constant) before all the incorrect beliefs vanish.

We note that many of the requirements that were needed for this theorem can be loosened (at the cost of obtaining weaker results). For example, we may get rid of the finiteness assumption on $\Theta$ (see \cite{cdc})  and we may assume the underlying graph is time-varying and directed (see \cite{acc1, acc2}).

\section{Simulations}  We now describe a quick simulation of our protocol for distributed optimization, using it to compute the median in linear time. The text and figures of this section are reproduced from \cite{linear}.

Specifically, each node in a network starts with a certain value $w_i$ and the nodes would like to 
compute a median of these values. Observe that the median of numbers $w_1, \ldots, w_n$ is a solution of the minimization problem
\[ \arg \min_{\theta} ~ \sum_{i=1}^n |\theta-w_i| \] and therefore can be computed in a distributed way using our distributed optimization method.

We implement the protocol of Eq. (\ref{optaccel}) on both the line graph and lollipop graph. In both cases, we  take $U=n$, i.e., we assume that every node knows the total number of nodes in the system. We suppose each node starts with $x_i(0)=w_i$. Moreover, for $i=1, \ldots, n/2$, we set $w_i$ to be the remainder of $i$ divided by $10$, and
$w_{n/2+i} = -w_i$.

We set the number of iterations $T$ as $T=4n$. In
Figure \ref{opt-simul}, we plot the number of nodes on the x-axis vs $(1/n) ||\widehat{y}(T)||_1$ (which is the average deviation from a correct answer since $0$ is a median) on the y-axis (here the hat-notation represents a running average).

As can be seen from the figures, choosing a linear number of iterations $T=4n$ clearly suffices to compute the median with good accuracy.

%%%%%%%%%%%%%%%%%%%%%%%%%%%%%%%%%%%%%%%%%%%%%%%%%%%%%%%%%%%%%%%%%%%%%%%%%%%%%%%%


\begin{thebibliography}{99}


\bibitem{degroot} M. H. DeGroot, ``Reaching a consensus,'' {\em Journal of the American Statistical
Association,} vol. 69, no. 345, pp. 118121, 1974






 \bibitem{JLM} A. Jadbabaie, J. Lin, and A. S. Morse, ``Coordination of groups of mobile autonomous
agents using nearest neighbor rules,''
{\em IEEE Transactions on Automatic
Control}, vol. 48, no. 6, pp. 988-1001, 2003.


 \bibitem{firstlearn} A. Jadbabaie, P. Molavi, A. Sandroni, and A. Tahbaz-Salehi, ``Non-bayesian social learning,'' {\em Games and Economic
Behavior}, vol. 76, no. 1, pp. 210–225, 2012.


\bibitem{ucsd} A. Lalitha, T. Javidi, and A. Sarwate, ``Social learning and distributed hypothesis testing,'' \url{https://arxiv.org/abs/1410.4307}










\bibitem{four} A. Nedic, A. Olshevsky, A. Ozdaglar, and J. N. Tsitsiklis, ``On distributed averaging
algorithms and quantization effects,'' {\em IEEE Transactions on Automatic
Control}, vol. 54, no. 11, pp. 2506-2517, 2009.


\bibitem{angas} A. Nedic, A. Ozdaglar, ``Distributed subgradient methods for multi-agent optimization,'' {\em IEEE Transactions on Automatic Control}, vol. 54, no. 1, pp. 48-61, 2009.


\bibitem{us} A. Nedic, A. Olshevsky, C. Uribe, ``Fast convergence rates for distributed non-Bayesian learning, '' \url{http://arxiv.org/abs/1508.05161}


\bibitem{acc1} A. Nedic, A. Olshevsky, C. Uribe, ``Nonasymptotic convergence rates for cooperative learning over time-varying directed graphs,'' {\em Proc. ACC 2015.}


\bibitem{cdc} A. Nedic, A. Olshevsky. C. Uribe, ``Distributed learning with infinitely many hypotheses,'' {\em Proc. CDC 2016}.


\bibitem{acc2} A. Nedic, A. Olshevsky, C. Uribe, ``Network independent rates in distributed learning,'' {\em Proc. ACC 2016.}




\bibitem{othesis} A. Olshevsky, ``Efficient Information Aggregation Strategies for Distributed Control and Signal Processing,'' Ph.D. thesis, MIT, Dept. of EECS, 2010.


\bibitem{linear} A. Olshevsky, ``Linear time average consensus on fixed graphs and implications for decentralized optimization and multi-agent control,'' preprint, \url{http://arxiv.org/abs/1411.4186}


%\bibitem{ram} S.S. Ram, A. Nedic, and V.V. Veeravalli %``Distributed stochastic subgradient projection %algorithms for convex optimization,'' {\em Journal of %Optimization Theory and Applications}, vol. 147, no. %3, pp. 516-545, 2010. 


\bibitem{rakh} S. Shahrampour, A. Rakhlin, and A. Jadbabaie, ``Distributed detection: Finite-time analysis and impact of network topology,''
\url{http://arxiv.org/abs/1409.8606}






\bibitem{tsi} J. N. Tsitsiklis, D. P. Bertsekas, and M. Athans, ``Distributed asynchronous deterministic
and stochastic gradient optimization algorithms,'' {\em IEEE Transactions
on Automatic Control,} vol 31, no. 9, pp. 803-812, Sep. 1986.


\end{thebibliography}
\end{document}